\def\version{May 21, 2013} 
\documentclass[12pt]{article}
\usepackage{amssymb}
\usepackage{hyperref}

\def\al{\alpha}

\def\claim{\par\bigskip\noindent{\bf Claim.} }
\def\de{\delta}
\def\ep{\epsilon}

\def\grp{{\Gamma}}
\def\la{\langle}
\def\ll{{\mathcal L}}
\def\om{\omega}

\def\poset{{\mathbb P}}
\def\proof{\par\noindent Proof\par\noindent}
\def\pr{\prime}
\def\qed{\par\noindent QED\par\bigskip}
\def\qq{{\mathbb Q}}
\def\ra{\rangle}
\def\res{\upharpoonright}
\def\rmand{\mbox{ and }}
\def\rr{{\mathbb R}}
\def\si{\sigma}
\def\sm{\setminus}
\def\span{\mbox{span}_\rr}
\def\st{\;:\;}
\def\su{\subseteq}
\def\summz{\sum_{n<\om}\zz}
\def\supp{{\mbox{supp}}}
\def\zspan{\mbox{span}_\zz}
\def\zz{{\mathbb Z}}

\newtheorem{theorem}{Theorem} 
\newtheorem{lemma}[theorem]{Lemma}
\newtheorem{question}[theorem]{Question}

\newtheorem{remark}[theorem]{Remark}

\begin{document}

\begin{center}
{\large
Countable subgroups of Euclidean space}
\end{center}

\begin{flushright}
Arnold W. Miller\\
April 2013 \\
revised \version  \\
\end{flushright}

\def\address{\begin{flushleft}
Arnold W. Miller \\
miller@math.wisc.edu \\
http://www.math.wisc.edu/$\sim$miller\\
University of Wisconsin-Madison \\
Department of Mathematics, Van Vleck Hall \\
480 Lincoln Drive \\
Madison, Wisconsin 53706-1388 \\
\end{flushleft}}

In his paper \cite{beros}, 
Konstantinos Beros proved a number of
results about compactly generated subgroups of Polish groups.
Such a group is $K_\si$, the countable union of compact sets.
He notes that the group of rationals under addition with the
discrete topology is an example of a Polish group which is
$K_\si$ (since it is countable) but not compactly generated
(since compact subsets are finite).

Beros showed that for any Polish group $G$, every $K_\si$
subgroup of $G$ is compactly generated iff every countable
subgroup of $G$ is compactly generated.  He showed that any
countable subgroup of $\zz^\om$ (infinite product of the integers)
is compactly generated and more generally,  for any Polish group
$G$, if every countable subgroup of $G$ is finitely generated,  then
every countable subgroup of $G^\om$ is compactly generated.

In unpublished work Beros asked the
question of whether
``finitely generated'' may be replaced by ``compactly generated''
in this last result.
He conjectured that the reals $\rr$ under addition
might be an example such that every countable
subgroup of $\rr$ is compactly generated but not
every countable subgroup of
$\rr^\om$ is compactly generated.
We prove (Theorem \ref{infdim}) that this is not true.
The general question remains open.

In the course of our proof we came up with some interesting
countable subgroups.   We show (Theorem \ref{plane}) that there
is a dense subgroup of the plane which meets every line in
a discrete set.  Furthermore for each $n$ there is a dense subgroup
of Euclidean space $\rr^n$ which meets every $n-1$ dimensional
subspace in a discrete set and a dense subgroup of $\rr^\om$ which
meets every finite dimensional subspace of $\rr^\om$
in a discrete set.

\begin{theorem}\label{reals}
Every countable subgroup $G$ of $\rr$ is compactly generated.
\end{theorem}
\proof
If $G$ has a smallest positive element, then this generates $G$.
Otherwise let $x_n\in G$ be positive and converge to zero.
Let $G=\{g_n\;:\;n<\om\}$.  For each $n$ choose $k_n\in \zz$
so that $|g_n-k_nx_n|\leq x_n$.
Let
$$C=\{0\}\cup \{x_n,\; g_n-k_nx_n \;\;:\;n<\om\}$$
then $C$ is a sequence converging to zero, so it is compact.
Also
$$g_n=(g_n-k_nx_n) + k_nx_n$$
so it generates $G$.
\qed

\begin{theorem}\label{ndim}
For $0<m<\om$ every countable subgroup $G$ of $\rr^m$ is compactly generated.
\end{theorem}
\proof
For any $\ep>0$ let
$$V_\ep=\span(\{u\in G\st ||u||<\ep\})$$
where here the span is taken respect to the field $\rr$.
Note that for $0<\ep_1<\ep_2$ that $V_{\ep_1}\su V_{\ep_2}$. And since
they are finite dimensional vector spaces 
there exists $\ep_0>0$ such that
$V_\ep=V_{\ep_0}$ whenever $0<\ep<\ep_0$.  Let $V_0=V_{\ep_0}$.
Let $\ep_0>\ep_1>\cdots$ descend to zero.  For each $n$ let
$B_n\su \{u\in G\;:\;||u||<\ep_n\}$ be a basis for $V_{\ep_n}=V_0$.
Then $|B_n|=\dim(V_0)$ and the sequence $(B_n)_{n<\om}$
``converges'' to the
zero vector $\vec{0}$.
Hence $\{\vec{0}\}\cup\bigcup_{n<\om}B_n$ is a compact subset of $G$.
Let
$$G_0=\zspan(\bigcup_{n<\om}B_n)$$
where $\zspan(X)$ is the set of
all linear combinations from $X$ with coefficients in $\zz$ or
equivalently the group generated by $X$.

\claim  $G_0$ is dense in $V_0$.

\proof
Let $|B_n|=k_0=\dim(V_0)$.
Suppose $v\in V_0$ and we are given $\ep_n$.
Then $v=\sum_{u\in B_n}\al_u u$
for some $\al_n\in\rr$.  Choose $n_u\in\zz$ with $n_u\leq \al_u<n_u+1$.
Let $w=\sum_{u\in B_n}n_uu_n$
Then
$$||v-w||=||\sum_{u\in B_n}(\al_u-n_u)u||\leq
\sum_{u\in B_n}|\al_u-n_u|\cdot ||u||
\leq |B_n|\ep_n=k_0\ep_n$$
Since $k_0\ep_n\to 0$ as $n\to\infty$ we are done.
\qed

Next we tackle the ``discrete part'' of $G$.

\claim There exists a finite $F\su G$ such that $G\su \zspan(F)+V_0$.

\proof
First note that for any $u,v\in G$, if $||u-v||<\ep_0$, then
$u-v\in V_0$ and hence $u+V_0=v+V_0$.
Now let $F_0\su G$ be finite so that
$$\span(F_0)=\span(G)$$
As we saw above for any $u\in G$ there exists
$v\in\zspan(F_0)$ such that 
 $$||u-v||\leq N_0=^{def}\sum_{u\in F_0} ||u||$$
If no such $F$ exists then we may construct an
infinite sequence $u_n\in G$ such that
  $$u_{n}\notin \zspan(F_0\cup\{u_i\st i< n\})+V_0
                 \rmand ||u_n||\leq N_0.$$
But by compactness of the
closed ball of radius $N_0$ there would have to
be $i<j$ with $||u_j-u_i||<\ep_0$ and
hence $u_j\in u_i+V_0$ which is a contradiction.
\qed

Now let
$$C=\{\vec{0}\}\cup F\cup \bigcup_{n<\om}B_n$$
and note that $C$ is a compact subset of $G$
and $\zspan(C)$ is dense in $G$.
To finish we can repeat the argument for $\rr$:

Let $G=\{g_n\st n<\om\}$.  Choose $u_n\in \zspan(C)$ with
$||g_n-u_n||<\frac1{2^n}$.  Then $C\cup\{g_n-u_n\st n<\om\}$
is compact and generates $G$.

\qed

\begin{question}
Do there exist Polish groups $\grp_1,\grp_2$ such that both
have the property that every countable subgroup is compactly generated
but the same is not true for $\grp_1\times\grp_2$?
\end{question}

Let $S$ be the real line with the Sorgenfrey topology,
i.e., half-open intervals $[x,y)$ are basic open sets.
Then $(S,+)$ is not a topological
group; addition is continuous but not subtraction.  However
countable subgroups are compactly generated since a descending sequence
with its limit point
is compact.  However $\{(x,-x):x\in S\}$ is discrete in $S\times S$, 
so the subgroup, so $\{(x,-x):x\in \qq\}$ is not compactly generated. 

\begin{theorem}\label{infdim}
Every countable subgroup $G\su \rr^\om$ is 
compactly generated.
\end{theorem}
\proof

Define $||x||_N=\sup\{|x(i)|: i<N\}$.  So this is just the sup-norm
but restricted to the first $N$-coordinates, i.e., $||x||_N=||x\res N||$.
An open neighborhood basis for $x$ is given
by the sets $\{y: ||y-x||_N<\ep\}$ for $N<\om$ and reals $\ep>0$.
Similarly we put 
$$G\res n=^{def}\{u\res n\::\; u\in\}.$$

For any $n$ by applying the proof of Theorem \ref{ndim} we may
find $\ep_n>0$ such that if 
$$V_n=^{def}\span(\{ u\res n\;:\;u\in G\rmand ||u||_n<\ep_n\})\su\rr^n$$
then for any $\de$ with $0<\de\leq\ep_n$
$$V_n=\span(\{ u\res n\;:\;u\in G\rmand ||u||_n<\de\}).$$
Without loss we may assume that $\ep_n<\frac1{2^n}$.
Take $B_n\su G$ so that
$$|B_n|=\dim(V_n)\leq n \rmand V_n=\span(\{u\res n\;:\;u\in B_n\}).$$
By the argument of the first claim of Theorem \ref{ndim},
we have that
for every $u$ with $u\res n\in V_n$
there exists $v\in \zspan(B_n)$
with 
$$||u-v||_n< n\cdot \ep_n<\frac{n}{2^n}$$
By the argument of the second claim there exists a 
finite $F_n\su G$ such that
$$G\res n\su \zspan (F_n\res n)+ V_n.$$
Hence for every $u\in G$ there exists 
$v\in\zspan(F_n\cup B_n)$ with $||u-v||_n<\frac{n}{2^n}$.

Apply this to each $u\in F_{n+1}$.  So we may choose
$\hat{u}\in \zspan(F_n\cup B_n)$ with $||u-\hat{u}||_n<\frac{n}{2^n}$.
Put 
$$H_{n+1}=^{def}\{(u-\hat{u})\;:\; u\in F_{n+1}\}$$
and put
$$C=\{\vec{0}\}\cup \bigcup_{n<\om} (B_n\cup H_n).$$
Then $C$ is compact, since for every $n$ if $u\in B_n$, then
we have that $||u||_n<\frac1{2^n}$ and if $v\in H_{n+1}$, then
$||v||_{n}<\frac{n}{2^{n}}$.
By induction we have that $F_n\su\zspan(C)$ for each $n$, since
if we assume $F_n\su\zspan(C)$, then for $u\in F_{n+1}$ we
have that $\hat{u}\in\zspan(C)$ and so $u=(u-\hat{u})+\hat{u}$
is in $\zspan(C)$.
It follows that $\zspan(C)$ is dense in $G$ since
$\bigcup_{n<\om} (B_n\cup F_n)$ is dense in $G$.
By the argument used in Theorems \ref{reals} and \ref{ndim}
we get that $G$ is compactly generated.
\qed

\begin{question}
For countable subgroups of $\rr^\om$ the generating compact set
can always be taken to be a convergent sequence.  We don't
know if this is more generally true. 
\end{question}

The motivation for the following example was a doomed
attempt to prove Theorem \ref{ndim} for the plane by looking
at one dimensional subspaces and considering multiple cases.

\begin{theorem}\label{plane}
There exists a dense subgroup $G\su \rr^2$ such that for
every line $L$ in the plane, $G\cap L$ is discrete.
\end{theorem}
\proof
Let $M$ be a countable transitive model of a large finite
fragment of ZFC.  Working in $M$ let $\poset$ be a countable
family of nonempty open subsets of the plane which is a basis
for the topology.   Forcing with $\poset$ produces a generic
point $p\in\rr^2$.

The following facts are well-known.
\begin{enumerate}
\item If $p$ is $\poset$-generic over $M$ then
 $p$ is not in any closed nowhere dense subset of $\rr^2$
 coded in $M$.
\item If $(p,q)$ are $\poset^2$-generic over $M$, then
 $p$ is $\poset$-generic over $M$ and
 $q$ is $\poset$-generic over $M[p]$.
\item If $(p,q)$ are $\poset^2$-generic over $M$, then
 $p+q$ is  $\poset$-generic over $M$.
\item If $p$ is $\poset$-generic over $M$ and $q\in M\cap\rr^2$,
 then $p+q$ is $\poset$-generic over $M$.
 And if $\al\in M\cap \rr$
 is nonzero, then $\al p$ is $\poset$-generic over $M$.
\end{enumerate}

Now let $\sum_{n<\om}\poset$ be the countable direct sum of
copies of $\poset$.  Let $(p_n\st n<\om)$ be 
$\sum_{n<\om}\poset$-generic over $M$ and let $G=\zspan(\{p_n\st n<\om\})$.

It is easy to see that the $p_i$ are algebraically independent
vectors over the field of rationals (or even the
field $M\cap\rr$) so for any $u\in G$
there is a unique
finite $F\su \om$ such that
let $u=\sum_{i\in F}n_ip_i$ and each $n_i\neq 0$ for $i\in F$.
Let $F=\supp(u)$ (the support of $u$).

Suppose $L$ is any line in the plane thru the
origin and let $$L^+=L\sm\{\vec{0}\}.$$
Note that $L$ is determined by any element of $L^+$.

\claim If $u,v\in L^+\cap G$, then $\supp(u)=\supp(v)$.

\proof
Let $\supp(u)=F_0$ and $\supp(v)=F_1$.
Suppose for contradiction that $F_1\sm F_0$ is nonempty.
Let $i\in F_1\sm F_0$ and put $F=F_1\sm\{i\}$.
Then
$$v=k_ip_i+\sum_{j\in F}k_jp_j.$$
It follows from the well-known facts above that
$p_i$ is $\poset$-generic over the model 
$N=^{def}M[(p_j :  j\in F_0\cup F)]$.
Also since $q=^{def}\sum_{j\in F}k_jp_j\in N$ and $k_i\in\zz$
we know that $v=k_ip_i+q$ is
$\poset$-generic over $N$.  On the other hand
the line $L$ is coded in $N$ since $u$ is
in this model.  But $L$ is a closed nowhere
dense subset of the plane so by genericity $v\notin L$.
\qed

Let $F_0$ be the common support of all $u\in L^+\cap G$.
Fix any $i_0\in F_0$.  Choose $u\in L^+\cap G$ such
and $u=\sum_{i\in F_0}n_ip_i$ and $n_{i_0}>0$ is the minimal
possible positive coefficient of $p_{i_0}$ and for any
point in $L^+\cap G$.

\claim $L\cap G=\{nu\st n\in \zz\}$.
\proof
Let $v=\sum_{i\in F_0}m_ip_i$ be in $L^+\cap G$.  Then
$m_{i_0}$ must be divisible by $n_{i_0}$ otherwise $n_{i_0}$
would not have been minimal.  Say $n\cdot n_{i_0}=m_{i_0}$
for some $n\in\zz$.  Then note that $v-n\cdot u$ is in $L\cap G$
and its support
is a subset of $F_0\sm\{i_0\}$.  Hence $v-n\cdot u=\vec{0}$
and we are done.
\qed

Now suppose that $L$ is a line not containing $\vec{0}$.  If
$L\cap G$ is nonempty take $u_0\in L\cap G$.  Let
$L_0=-u_0+L$.  Since $L_0$ is line thru the origin there
exists $u_1\in G$ such that
$G\cap L_0=\{nu_1\st n\in\zz\}$.  Then
$$G\cap L=\{u_0+nu_1\st n\in\zz\}.$$

\qed

\begin{remark}.

\begin{enumerate}
\item
Note that dense-in-itself is not the same as dense, for example,
$\zz\times\qq$ is dense-in-itself but not dense in $\rr^2$.
\item  Finitely generated does not imply discrete.  The additive subgroup
of $\rr$ generated by $\{1,\sqrt{2}\}$ is dense.
\item  A circle $C$ is closed nowhere dense subset of the plane and
determined by any three points, hence $C\cap G\su \zspan(F)$
for some finite $F\su G$.
Also the graph $Q$ of a polynomial of degree $n$
is determined by any $n+1$ points so there will be a finite $F\su G$
with $G\cap Q\su \zspan(F)$.
\item  One can get a group $G\su\rr^2$ which is continuum-dense
and meets every line in a discrete set.  To see this
construct a sequence of perfect
sets $P_n$ such that every nonempty open set contains
one of them and such
that for every finite $\{p_i: i<n\}\su\cup_nP_n$
of distinct points the tuple
$(p_i: i<n)$ is $\poset^n$-generic over $M$.
Then $G=\zspan(\cup_n P_n)$ meets every line discretely.

\end{enumerate}
\end{remark}

\begin{theorem}\label{dense}
There exists a dense subgroup $G\su\rr^N$ such that
$G\cap V$ is discrete for
every $V\su\rr^N$ a vector space over $\rr$ of dimension
$N-1$.
\end{theorem}
\proof

In the model $M$ we let $\poset$ be the nonempty open subsets of
$\rr^N$.  Take $(p_i:i<\om)$ be $\sum_{n<\om}\poset$-generic over
$M$ and put $G=\zspan(\{p_n:n<\om\})$.

Let $\summz$ be the sequences  $x\in\zz^\om$
such that $x(n)=0$ except
for finitely many $n$.
Let $h:\summz\to \rr^N$ be the homomorphism defined
by
$$h(x)=\sum_{n<\om}x(n)p_n$$
and note that this is a finite sum, the coefficients are integers
and the $p_n$ are generic vectors in $\rr^N$.

Clearly the generic vectors are linearly independent over
the rationals so the map $h$ is a group isomorphism
$h:\summz\to G$.

Let $C\su \summz$ be $h^{-1}(G\cap V)$.
As in the proof of Theorem 4.2 Beros \cite{beros}, let
$$C_n=\{x\in C\st \forall i< n \;\; x(i)= 0\}.$$
If $C_{n+1}$ is a proper subgroup of $C_n$
choose $k_n\in\zz$ non-zero which
divides every $y(n)$ for $y\in C_n$.  Otherwise,
if $C_{n+1}=C_n$ put $k_n=0$.
For each $n$ choose $x_n\in C_n$ with $x_n(n)=k_n$.

\claim If $k_n\neq 0$, then $h(x_n)\notin \span(\{h(x_i):i>n\})$.

\proof
Suppose for contradiction that for some $m$
$$h(x_n)\in\span(\{h(x_{n+i}):i=1\ldots m\}).$$
Let $V_0=\span(\{h(x_{n+i}):i=1\ldots m\})$ and note
that $V_0\su V$ and so $\dim(V_0)\leq \dim(V)<N$.
Hence $V_0$ is closed and nowhere dense in $\rr^N$.

Now $\la h(x_{n+i}):i=1\ldots m \ra\in M[\la p_i : i>n\ra ]$
and hence $V_0$ is a closed subset of $\rr^N$ coded
in $M[\la p_i : i>n\ra ]$.  But
$$h(x_n)=k_np_n+\sum_{i>n}x_n(i)p_i$$
and so as we saw in the last proof
$h(x_n)$ is $\poset$-generic over $M[\la p_i : i>n\ra]$.
But this means that it avoids
the closed nowhere dense set $V_0$ coded in $M[\la p_i : i>n\ra]$.
Contradiction.
\qed

\claim $|\{n:k_n\neq 0\}|\leq \dim(V)$.
\proof
Suppose not and let $Q\su \{n:k_n\neq 0\}$ be finite with
$|Q|>\dim(V)$.
Then for each $n\in Q$ we have that
$$h(x_n)\notin\span(\{h(x_m)\st m\in Q \rmand m>n\}). $$
It follows that $\{h(x_n):n\in Q\}$ is a linearly independent
subset of $V$ of size greater than $\dim(V)$, which
is a contradiction.
\qed

Hence $Q=^{def}\{n: k_n\neq 0\}$ is finite.  This means
that
$$C=\zspan(\{x_n:n\in Q\})$$
by a process resembling Gaussian elimination.  Knowing that
$Q$ is finite shows that the algorithm terminates.

It follows that
$V\cap G=\zspan(\{h(x_n):n\in Q\})$
and $\{h(x_n):n\in Q\}$ is a basis for $\span(V\cap G)$.
Extend $\{h(x_n):n\in Q\}$ to a basis $B$ for $\rr^N$.
Take a bijective linear transformation
$j:\rr^N\to \rr^N$
determined by mapping the standard basis of $\rr^N$ to $B$.
It follows that $j(\zz^{|Q|}\times\{\vec{0}\})=V\cap G$
and so $V\cap G$ is discrete in $\rr^N$.
\qed

\begin{remark}
For $\rr^\om$ using this technique 
we get a dense generically generated group
$G\su \rr^\om$ such that $V\cap G$ is discrete for every finite dimensional
subspace $V\su\rr^\om$.   
\end{remark}

Although we are enamored with the generic sets argument
of Theorem \ref{dense}, it also possible 
to give an elementary inductive construction.  We haven't
verified it in general but we did check the proof for the
plane.  As long as we are at it, we might as well construct a
subgroup of $\rr^2$ which is maximal with respect to
linear discreteness,
although this will require a subgroup of cardinality the
continuum.

\begin{theorem}
Let $\ll$ be the family of all lines in the plane containing
the origin.
There exists a subgroup $G\su\rr^2$ such that
$G\cap L$ is an infinite discrete set for
every $L\in\ll$.
\end{theorem}
\proof
Define $G\su\rr^2$ is linearly discrete iff $G\cap L$ is discrete
for every $L\in\ll$.

\noindent Define $\ll_G=\{L\in \ll\;:\; |G\cap L|\geq 2\}$.

\begin{lemma}
For any linearly discrete subgroup $G\su\rr^2$ with $|G|<|\rr|$ and
any $L_0\in \ll\sm\ll_G$ there exists $p\in L_0\sm\{\vec{0}\}$ such
that $G+\zz p$ (the group generated by $G$ and $p$)
is linearly discrete and $(G+\zz p)\cap L=G\cap L$
for every $L\in\ll_G$.
\end{lemma}
\proof
Suppose $np+q\in L$ where $n\in\zz$ is nonzero,
$L\in\ll_G$, and $q\in G$.  Then $p+\frac{q}n\in L$ since $L$ contains
the origin.  But the lines $L_0$ and $L-\frac{q}n$ meet in at most one
point.  It follows that at most $|G|$ many points of $L_0$ are
ruled out so we may choose $p\in L_0\sm\{\vec{0}\}$ such
$(G+\zz p)\cap L=G\cap L$
for every $L\in\ll_G$.

This implies that $G+\zz p$ is linearly discrete by the following
argument.
We only need to worry about new lines $L\in\ll$ which are
not in $\ll_G$.  Suppose $np+q\in L$ where $L\in\ll\sm\ll_G$,
$n\in\zz$ nonzero, and $q\in G$.
Choose such a point on $L$ with minimal $n>0$.
We claim that $L\cap (G+p\zz)=\{knp+kq\:;\; k\in\zz\}$.

Then for any $k\in \zz$
we have that $knp+kq\in L$.
Any point
on $L\cap (G+p\zz)$ has the form 
Let $mp+q^\pr$ be an arbitrary point of $p\zz+G$.  
If it is in $L$ it must be that $n$ divides $m$ otherwise 
$0<m-kn<n$ for some
$k\in\zz$ and 
$$(mp+q^\pr)-(knp+kq)=(m-kn)p+(q^\pr-kq)\in L$$ 
shows that $n$ was not minimal.  Hence
for some $k\in\zz$ we have that $kn=m$.  It must be
that $kq=q^\pr$ since otherwise 
$q^\pr-kq\in L$ is a nontrivial element of $G$ in $L$, contradicting
$L\notin \ll_G$.
\qed

Using the Lemma it is an easy inductive construction of an
and increasing sequence of subgroups $G_\al\su\rr^2$ for $\al<|\rr|$ 
such that $|G_\al|=|\al+\om|$ is linearly discrete and
for each line $L\in\ll$ there is an $\al$ with 
$|G_\al\cap L|\geq 2$.
\qed

\begin{question}
For a line $L$ in the plane not containing the origin,
if $G$ meets $L$ at all, then taking $u\in G\cap L$ and considering
$L-u\in\ll$ shows us that $G$ meets $L$ in an infinite discrete
set.  I don't know if we can construct such a $G$ which meets every
line.  I don't even know how to meet all of the horizontal
lines.
\end{question}

Jan Pachl points out that every locally compact compactly generated
Abelian group is a product of a compact Abelian group and a finite number of
copies of $\rr$ and $\zz$ (see Hewitt and Ross \cite{hewitt} 9.8).
Also,
every locally compact Abelian group is a product of $\rr^k$ and 
a locally compact
subgroup that has a compact open subgroup (see Hofmann and Morris
\cite{hofmann} 7.57).  Every
$K_\si$ topological group is a closed subgroup of a 
compactly generated group (see Pestov \cite{pestov}).

\bigskip

\address

\end{document}